>


\baselineskip=14pt
\parskip=10pt

\magnification=\magstephalf

\def\C{{\cal C}}

\def\1{{\overline{1}}}
\def\2{{\overline{2}}}
\parindent=0pt
\overfullrule=0in

\def\frac#1#2{{#1 \over #2}}
\centerline
{\bf
Some Deep and Original Questions about the ``critical exponents'' of Generalized Ballot Sequences
 }

\bigskip
\centerline
{\it Shalosh B. EKHAD and Doron ZEILBERGER}

\bigskip

{\bf Abstract}: We numerically estimate the critical exponents of certain enumeration sequences
that naturally  generalize the famous Catalan and super-Catalan sequences, and raise
deep and original questions about their exact values, and whether they
are rational numbers.

{\bf Added  April 7, 2021}: Michael Wallner answered the questions raised in this article (in the 3D case).
Once he writes it up, we will reference it.

{\bf Added  May 27, 2021}: Michael Wallner just posted his beautiful article

https://arxiv.org/abs/2105.12155

The promised donation to the OEIS in his honor has been made. Wallner also pointed out some minor typos. For the sake of historical
accuracy, we did not correct them, and this second version of our paper is the same as the original one, after these remarks.

{\bf Preface}

Some of the deepest, and most original, questions in enumeration are those of finding the so-called
{\it critical exponents} of hard-to-compute sequences, most famously sequences enumerating self-avoiding
walks. Physicists believe that the critical exponent  is even more interesting than the so-called connective constant
(since it is believed to be `universal' and tells you `how water boils', as opposed to `what is the boiling temperature?', the
latter depending on many contingent factors).

An important kind of sequences that occur a lot in enumerative combinatorics is the class of  {\bf P-recursive sequences}
for which it is easy to derive exactly both connective constants (see below, the limit of the ratio of consecutive terms)
and {\bf critical exponents} (see below). 
In this deep and original article, we come up with many interesting
and natural enumerative sequences, generalizing the Catalan numbers and their higher-dimensional counterparts
(ballot numbers), that most probably are {\bf not} $P$-recursive, yet seem to have nice
asymptotics, similar to those of $P$-recursive sequences (and self-avoiding walks). We  computed
quite a few of them and estimated their critical exponents, but we can't even conjecture exact values,
and whether or not they are rational. We are pledging donations to the OEIS in honor of the
first identifier (and prover) of these critical exponents.

{\bf Maple Package}

This article is accompanied by the Maple package {\tt Capone.txt}, available from,

{\tt https://sites.math.rutgers.edu/\~{}zeilberg/tokhniot/Capone.txt} \quad.

The front of this article

{\tt https://sites.math.rutgers.edu/\~{}zeilberg/mamarim/mamarimhtml/capone.html} \quad, 

contains numerous input and output files, briefly described in this article.

{\bf Ballot Sequences}

One of the almost $\aleph_0$ combinatorial objects counted by the {\it Catalan  numbers}, $\C_2(n):=\frac{(2n)!}{n!(n+1)!}$,
are {\it ballot sequences}, that count the number of ways of counting votes of two candidates that
each got $n$ votes, in such a way that one of them was  never behind in the partial count, as first proved
by Bertrand [B], and given the {\it proof from the book} by Andr\'e [A]. It is one of the {\it fundamental objects}
listed in Stanley's fascinating book ([St], Theorem 1.5.1 (iv)).

But what if  there are $k$ candidates? What is the number of ways of counting the votes of $k$ candidates, each
of whom got at the end $n$ votes, in such a way that the second candidate was never ahead of the first candidate,
the third was never ahead of the second one, and in general, for $i=1 \dots k-1$, Candidate $i+1$ was never ahead
of Candidate $i$? Equivalently:

{\it ``How many ways are there of walking in the $k$-dimensional Manhattan lattice, with unit steps in each direction,
from the origin $(0, \dots, 0)$ to the point $(n, \dots, n)$ always staying in the region
$$
\{ (x_1, \dots, x_k) \, | \, x_1 \geq x_2 \geq \dots \geq x_{k-1} \geq x_k \}?'' \quad
$$
}

This number is famously
$$
C_k(n):=\frac{(n\,k)! \, 0! \, 1! \, \cdots \, (k-1)!}{(n+k-1)!\, (n+k-2)!\, \dots\, n!} \quad.
$$

Note that when $k=2$, one gets the Catalan numbers.

The formula for $C_k(n)$ follows immediately from the so-called {\it Young-Frobenius} formula, equivalent
to the {\it hook-length formula} (see [Kn], p. 58, Eq. (34)). The latter formula tells you that the
number of ways of walking from the origin to $(a_1, \dots, a_k)$ (when $a_1 \geq a_2 \geq \dots \geq a_k$)
{\bf always} staying in $\{ (x_1, \dots, x_k) \, | \, x_1 \geq x_2 \geq \dots \geq  x_{k-1} \geq x_k \} \quad$ is
$$
\Delta(a_1+k-1, a_2+k-2, \dots, a_k)\, \cdot \, \frac{(a_1 + \dots + a_k)!}{(a_1+k-1)!\, (a_2+k-2)! \, \cdots \, a_k!} \quad ,
$$
where $\Delta(x_1, \dots , x_k)$ is the {\it discriminant} function 
$$
\Delta(x_1, \dots , x_k) \, := \, \prod_{1 \leq i < j \leq k} \, (x_i -x_j) \quad .
$$

Substituting $a_1=n, a_2=n, \dots, a_k=n$ immediately gives the above expression for $C_k(n)$.

What's nice about the Catalan numbers and their higher dimensional analogs is that they
belong to a very important class of sequences, called {\bf   $P$-recursive sequences}. 

A sequence $x(n)$ is $P$-recursive if there exists a {\bf finite} positive integer $L$, called its {\bf order}, and
{\bf polynomials} $p_i(n)$, for $i=0,1, \dots, L$, such that
$$
p_0(n) x(n) \, + \,p_1(n) x(n+1) \, + \, \dots \,+ \,  p_L(n) x(n+L) \, = \, 0 \quad.
$$

It so happens that for the Catalan sequence $C(n)=C_2(n)$, and even the super-Catalan sequences $C_k(n)$,
the order $L$ is $1$. Such sequences are called {\bf hypergeometric sequences}. But from a computational
point of view, for {\it any} such sequence, once you have found a recurrence satisfied by it, you can compute the first $10000$ (or whatever) terms very fast.

It is easy to see, using elementary linear algebra, that if an integer sequence is $P$-recursive, then
the coefficients of the polynomials $p_i(n)$ above are all integers. This observation will
be important later on.

It is easy to see ([KaP]) that a sequence $x(n)$ is $P$-recursive if and only if its {\bf generating function}, alias
{\bf z-transform},
$$
X(z) \, := \, \sum_{n=0}^{\infty} x(n) \, z^n \quad,
$$
satisfies a {\bf linear differential equation with polynomial coefficients}. The sequence $\{x(n)\}$ is called
{\bf algebraic} if its generating function $X(z)$ satisfies an {\bf algebraic equation},
i.e., for some positive integer $M$, called the {\bf degree}, we have
$$
q_0(z) \, + \, q_1(z) X(z) + q_2(z) X(z)^2 \, + \dots \,+ q_M(z) X(z)^M \, = \, 0 \quad,
$$
for some {\bf polynomials} of $z$, $q_0(z), ..., q_M(z)$. Once again it is easy to see that for integer sequences that are algebraic,
the coefficients of the $q_i(z)$ can be taken to be integers.

If the degree of the defining equation, $M$, is $1$, then the generating function is {\bf rational}
and the sequence belongs to the important subfamily of {\bf $C$-finite sequences} (see [KaP], chapter 4).

The Catalan sequence $C(n)$ is  algebraic since it famously satisfies the equation
$$
z\,X(z)^2 \,- \,X(z) \,+ \, 1 \, = \, 0 \quad .
$$

Perhaps surprisingly, the higher-dimensional counterparts, $C_k(n)$, for $k \geq 3$, are {\bf not} algebraic.

A classical theorem (that may be due to Comtet, see [KaP]) states that if a sequence has an algebraic
generating function, then it is $P$-recursive, but of course, not vice-versa, witness the fact
that $C_k(n)$ are {\bf not} algebraic for $k \geq 3$.

Using standard techniques it is easy to prove the following fact.

{\bf Important Fact}: Every $P$-recursive integer sequences, $x(n)$, has nice asymptotic expressions and
there exist algorithms for finding it. For most `natural' sequences the asymptotics is of the form
$$
x(n) \asymp C \mu^n n^{\theta} \quad,
$$
where $C$ is a constant, and, borrowing the terminology of mathematical physics, $\mu$ is called the {\bf connective constant}
and $\theta$ is called the {\bf critical exponent}.

{\bf The Really Deep Sequences: Sequences enumerating Self-Avoiding Walks}

Tony Guttmann, Iwan Jensen, and their collaborators used sophisticated (alas, still {\it exponential-time}) algorithms
to crank out many terms of the sequences of self-avoiding walks in the 2-dimensional, triangular, square, and honeycomb
lattices. Just search the OEIS for "Self-avoiding walks``. The only rigorously proved connective constant
is the one for the number of self-avoiding walks on the honeycomb lattice, proved to be $\sqrt{2+\sqrt{2}}$ in [DS]
confirming a previous non-rigorous, but very convincing, `physical derivation' by B. Nienhuis [N].

The {\bf exact value} of the critical exponent for {\bf all} (hence the {\it universality}) 
two-dimensional lattices (the triangular, cubic, and honeycomb) is {\bf conjectured} to be $\frac{11}{32}$
(the physicists, for their own reasons,  add $1$, so for them it is $\frac{43}{32}$), but this
is {\bf wide open}, and worthy of a Fields medal.

{\bf Generalized Ballot Sequences in 2 dimensions}

The general question is to enumerate the number of walks, with positive unit steps, in the two-dimensional square lattice
from the origin to the point $(an, bn)$ always staying in the region  $b\,x_1\,-\, a\,x_2 \geq 0$. Of course $a$ and
$b$ are relatively prime. The special cases $a=1$, general $b$ 
(and $b=1$, general $a$)  give the {\bf Fuss-Catalan numbers}, and were
considered by R.C. Lyness (of ``Lyness cycle'' fame, one of the greatest high school teachers of all time) [L].
It is known (e.g. [BKK] [Ek] [TZ]) that the generating function is {\bf always} algebraic, and there
are efficient Maple implementations in [EK] and [TZ] for finding them. The connective constant is
easily seen to be
$$
\frac{(a+b)^{a+b}}{a^a \, b^b} \quad.
$$
It turns out that the critical exponent is {\bf always}  {\bf exactly} $-\frac{3}{2}$, and it is not hard to prove this using standard
techniques.

To get the first $100$ terms, as well as precise asymptotic {\it estimates} for {\bf all} such sequences for $1 \leq a <b \leq 10$ and
relatively prime $a$, and $b$, see the output file

{\tt http://www.math.rutgers.edu/\~{}zeilberg/tokhniot/oCapone1.txt} \quad  .

Some of these sequences (for small $a$ and $b$) are in the OEIS.

The estimated critical exponents, using our {\bf numerical method}, agree to (at least) ten decimal digits with
the rigorously derived value of $-\frac{3}{2}$, giving us confidence at the estimates for the higher dimensional 
generalized ballot sequences, to be discussed next, for which, at present, there is {\bf no} rigorous way (as far as we know).

{\bf Generalized Ballot Sequences in 3 dimensions}

Things are starting to get interesting in three dimensions. To our surprise, {\bf none} of these sequences
were in the OEIS (viewed April 4, 2021).

The output file

{\tt http://www.math.rutgers.edu/\~{}zeilberg/tokhniot/oCapone2.txt} \quad  ,

contains enumerating sequences for the number of ways of walking from $(0,0,0)$ to $(a\,n \,,\, b\, n \, , \, c \, n)$ 
such that if $M=lcm(a,b,c)$, it stays in the region  $(M/a)\, x \geq (M/b)\, y \geq (M/c)\,z \geq 0$ for
all $1 \leq c \leq b \leq a \leq 4$ such that $gcd(a,b,c)=1$.

We believe that except for the classical case $a=b=c=1$, these sequences are {\bf not} P-recursive, and hence
we have no clue how to derive, rigorously, the critical exponents. We don't even know
whether they are rational numbers, all we can do is crank out many terms and use {\bf numerics}.

Here are numerical estimates for $(a,b,c)$, as above, taken from the end of the above output file.

$\bullet$ (1,1,1) [the classical (proved!) case]: $-4$ \quad .

$\bullet$ (2,1,1):   $-3.7312$ \quad .

$\bullet$ (2,2,1):   $-4.2884$ \quad .

$\bullet$ (3,1,1):   $-3.5976$ \quad .

$\bullet$ (3,2,1):   $-4.055$ \quad .

$\bullet$ (3,2,2):   $ -3.8375$ \quad .

$\bullet$ (3,3,1):   $-4.4455$ \quad .

$\bullet$ (3,3,2):   $-4.1695$ \quad .

$\bullet$ (4,1,1):   $-3.515$ \quad .

$\bullet$ (4,2,1):   $-3.9091$ \quad .

$\bullet$ (4,3,1):   $-4.2454$ \quad .

$\bullet$ (4,3,2):   $-4.0237$ \quad .

$\bullet$ (4,3,3):   $-3.8834$ \quad .

$\bullet$ (4,4,1):   $-4.5453$ \quad .

$\bullet$ (4,4,3):   $-4.12019$ \quad .

{\bf Disclaimer}: We have no rigorous error bars for the above values, so they may be only right to fewer decimal places.
We used our own (admittedly crude) {\it home-made} asymptotics, implemented in procedure 

{\tt MyAsyM(L,mu,n,k)} \quad ,

that inputs a sequence of positive numbers (integers in our case) believed to have an asymptotic behavior of the form 
(where $\mu$, luckily, is known beforehand)
$$
C \, \mu^n n^\theta (1+ \frac{c_1}{n}+ \frac{c_2}{n^2} + \dots )
$$
and pretending that the asymptotic expansion only goes as far as the term $\frac{c_k}{n^k}$, takes the $\log$ of the last few terms of the sequence and
solves for the unknowns $\log C$, $\theta$ and $c_1,c_2, \dots, c_k$. Of course, our main interest is in $\theta$.

Even though our algorithms are {\it polynomial time} (essentially dynamical programming), computing the first $N$ terms is $O(N^3)$ (for the cubic lattice case).
To get a more reliable estimate for the case $(2,1,1)$ we computed $400$ terms in the output file

{\tt http://www.math.rutgers.edu/\~{}zeilberg/tokhniot/oCapone4.txt} \quad  ,

that yields the estimate $-3.731220575$ for the critical exponent for that sequence. We would love to know whether the
exact value is a rational number, and are offering a donation of 100 dollars to the OEIS in honor of the first (proved) exact answer.

{\bf Generalized Ballot Sequences in 4 dimensions}

Things start to get slow for higher dimensions, but for a few cases see this output file:

{\tt http://www.math.rutgers.edu/\~{}zeilberg/tokhniot/oCapone3.txt} \quad  .

{\bf Sequences enumerating the total number of $n$-step walks}

The total number of  $n$-step walks in the $k$-dimensional hyper-cubic lattice
that stay in the region $x_1 \geq x_2 \geq \dots \geq x_k \geq 0$ for, $2 \leq k \leq 12$,  are all in the OEIS.
These also count  $n$-celled Standard Young tableaux with at most $k$ rows, and
thanks to the Robinson-Schenstead correspondence ([Kn]), the number of $12 \dots (k+1)$-avoiding {\it involutions} of size $n$.

$\bullet$  $k=2$,  A1405 ({\tt https://oeis.org/A001405}). Easily seen to be given by ${ {n} \choose {[ n/2]}}$.

$\bullet$  $k=3$,  A1006 ({\tt https://oeis.org/A001006}). First discovered, and proved,  by Amitai Regev [R],
to be given by the almost-as-famous Motzkin numbers.

Higer dimensional cases were treated in [BFK].

$\bullet$  $k=4$: A5817 ({\tt https://oeis.org/A005817}, see the references there).

$\bullet$  $k=5$: A49401, ({\tt https://oeis.org/A049401}).

$\bullet$  $k=6$: A7579, ({\tt  https://oeis.org/A007579}).

$\bullet$  $k=7$: A7578, ({\tt  https://oeis.org/A007578}).

$\bullet$  $k=8$: A7580, ({\tt  https://oeis.org/A007580}).

$\bullet$  $k=9$: A212915, ({\tt  https://oeis.org/A212915}).

$\bullet$  $k=10$: A212916, ({\tt  https://oeis.org/A212916}).

$\bullet$  $k=11$: A229053, ({\tt  https://oeis.org/A229053}).

$\bullet$  $k=12$: A229068, ({\tt  https://oeis.org/A229068}).

$\bullet$ $k=13$ was not in the OEIS (viewed April 4, 2021). The first 16 terms are

$$
1, 2, 4, 10, 26, 76, 232, 764, 2620, 9496, 35696, 140152, 568504, 2390479, 10349521, 46206511
$$

It is well-known that these sequences are  $P$-recursive for all $k$, so once a recurrence is found it is
easy to get many terms, as well as precise asymptotics.

But what about other regions?

Even for the two-dimensional case, except for the sequence enumerating the number of $n$-step walks
in the 2D  square lattice  that stay in the region $\{ (x_1, x_2) \, | \, x_1 \geq 2 x_2 \, \}$, that
is OEIS sequence  A126042 ({\tt http://oeis.org/A126042}), that is there for a different reason
(and it would be an interesting exercise to prove that they are indeed the same), {\bf none} of the
other cases, even for the 2D case, seem to be there (yet).

Once again these (in the 2D case) are all  algebraic, and henec $P$-recursive. We conjecture that for higher dimensions
they are not $P$-recursive in general.

For the $2$-dimensional case, see the output file:

{\tt https://sites.math.rutgers.edu/\~{}zeilberg/tokhniot/oCapone5.txt} \quad  .

For the $3$-dimensional case, see the output file:

{\tt https://sites.math.rutgers.edu/\~{}zeilberg/tokhniot/oCapone6.txt} \quad .

Here the critical exponents (separated according to $n$ mod $a+b+c$) are (conjecturally!) {\it nice}:  either $-\frac{1}{2}$ or $-1$.

{\bf Conclusion: Why is this Research both Deep and Original?}

Determining the (proved!) {\bf exact values} of critical exponents of enumeration sequences that do not seem to be $P$-recursive is very challenging.
In the case of {\it self-avoiding walks}, it may get you a Fields medal. Since this is out of reach at present, it is fun and interesting to experiment
with other kinds of walks, that are just minor tweaks of the classical ballot sequences, yet seem much harder, and we suspect
that for dimensions three and higher are {\bf not} $P$-recursive. But perhaps they belong to another class 
for which it would be possible to find the exact asymptotics? What's nice about the sequences considered
in the present paper is that the algorithms for generating many terms of the sequences are {\bf polynomial time}, hence with larger computers, one
should be able to get many more terms. and derive more precise estimates.

Another kind of walks was considered in [KaZ], where
the region was the classical one $x_1 \geq x_2 \geq \dots \geq x_k$ but certain runs were forbidden. Surprisingly
there the (conjectured!)  critical exponents turned out, often, to be very simple. 

If nothing else, we found lots of new sequences! Hopefully some of them will be entered in the OEIS by
kind readers.

{\bf References}

[A] D. Andr\'e, {\it Solution directe du probl\`eme r\'esolu par M. Bertrand}, Comptes Rendus Acad. Sci. Paris {\bf 105} (1887), 436-437.

[BFK] F. Bergeron, L. Favreau and D. Krob, 
{\it Conjectures on the enumeration of tableaux of bounded height}, Discrete Math {\bf 139} (1995), 463-468.

[B] J. Bertrand, {\it Solution d'une probl\`eme}, Comptes Rendus Acad. Sci. Paris {\bf 105} (1887), 369.

[BKK] Cyril Banderier, Christian Krattenthaler, Alan Krinik, Dmitry Kruchinin, Vladimir
Kruchinin, David Tuan Nguyen, and Michael Wallner, {\it Explicit formulas for enumeration
of lattice paths: basketball and the kernel method}, 2017. \hfill\break
{\tt https://arxiv.org/abs/1609.06473v2}

[DS] H. Domini-Chupin and S. Smirnov, {\it The connective constant of the honeycomb lattice is $\sqrt{2+\sqrt{2}}$},
{\tt https://arxiv.org/abs/1007.0575}. \hfill\break
Also in: Anal. of Math. (2) {\bf 175} (2012), 1653-1665.

[Ek] Bryan Ek, {\it Lattice Walk Enumeration}, {\tt https://arxiv.org/abs/1803.10920}

[KaP] Manuel Kauers  and Peter Paule, {\it ``The Concrete Tetrahedron''}, Springer, 2011.

[KaZ] Manuel Kauers and Doron Zeilberger, {\it Counting Standard Young Tableaux with restricted runs},
Personal J. of Shalosh B. Ekhad and Doron Zeilberger. \hfill \break
{\tt https://sites.math.rutgers.edu/\~{}zeilberg/mamarim/mamarimhtml/cyt.html} \quad .

[Kn] Donald E. Knuth, {\it ``The Art of Computer Programming''}, Second Edition, Addison-Wesley, 1998.

[L] R.C. Lyness, {\it Al Capone and the Death Ray}, The Mathematics Gazette {\bf 25} (1941), 283-287. \hfill\break
[available from JSTOR]

[N] B. Nienhuis, {\it Exact critical point and critical exponents of the $O(n)$ models in two dimensions},
Phys. Rev. Letters {\bf 49} (1982), 1062-1065.

[R] Amitai Regev, {\it Asymptotic values for degrees associated with strips of Young diagrams},
Adv. in Math. {\bf 41} (1981),  115-136.

[Sl] Neil Sloane, {\it The On-Line Encyclopedia of Integer Sequences}, {\tt  https://oeis.org} \quad .

[St] Richard P. Stanley, {\it ``Catatan numbers''}, Cambridge University Press, 2015.

[TZ] Thotsaporn ``Aek'' Thanatipanonda, Doron Zeilberger, 
{\it A Multi-Computational Exploration of Some Games of Pure Chance}, {\tt https://arxiv.org/abs/1909.11546}. \hfill\break
Also in: J. of Symbolic Computation {\bf 104} (2021), 38-68.

[WimZ] Jet Wimp and Doron Zeilberger, 
{\it Resurrecting the asymptotics of linear recurrences}, 
J. Math. Anal. Appl. {\bf 111}, 162-177 (1985).


\bigskip
\hrule
\bigskip
Shalosh B. Ekhad and Doron Zeilberger, Department of Mathematics, Rutgers University (New Brunswick), Hill Center-Busch Campus, 110 Frelinghuysen
Rd., Piscataway, NJ 08854-8019, USA. \hfill\break
Email: {\tt [ShaloshBEkhad, DoronZeil] at gmail dot com}   \quad .

{\bf Exclusively published in the Personal Journal of Shalosh B. Ekhad and Doron Zeilberger and arxiv.org}

Original version, written: {\bf April 4, 2021}.

This version, written: {\bf May 27, 2021}.

\end